\documentclass[reqno]{amsart}
\begin{document}

\title[]{Remark on Cavalieri's Quadrature Formula}

\author{David M. Bradley}
\address{Department of Mathematics and Statistics\\
         University of Maine\\
         5752 Neville Hall\\
         Orono, Maine 04469--5752\\
         U.S.A.}
\email{bradley@math.umaine.edu, dbradley@member.ams.org}

\date{Submitted
November 2, 2002.}

%\subjclass{Primary: 44A20; Secondary: 65D30, 65R10}

%\keywords{Integral Transforms, Simpson's Rule, Estimating
%Derivatives, Differentiating under the Integral Sign.}

\maketitle

%\section{Introduction}
With regard to Cavalieri's quadrature formula~\cite{Wild} for
$0<a<b$ and integer $n$:
\begin{equation}
   \int_a^b x^n \,dx = \begin{cases} \frac{b^{n+1}-a^{n+1}}{n+1},
   \quad &\mbox{if $n\ne -1$},\\
    \log\frac{b}{a}, \quad &\mbox{if $n=-1$},\end{cases}
\label{cav}
\end{equation}
there is an amusing sense in which the exceptional case contains
the others.  For example, if $|t|<1$ then
\[
   \sum_{n=0}^\infty t^n\int_0^1 x^n \,dx
   = \int_0^1 \sum_{n=0}^\infty (tx)^n
   = \int_0^1 \frac{dx}{1-tx}
   = t^{-1}\log(1-t)^{-1}
   = \sum_{n=0}^\infty\frac{t^n}{n+1}.
\]
More generally, if $0<a<b$ and $|t|<\min(a,1/b)$ then
\[
   \sum_{n=0}^\infty t^n\int_a^b x^n\,dx =
   t^{-1}\log\bigg(\frac{1-at}{1-bt}\bigg), \qquad
   \sum_{n=1}^\infty t^n\int_a^b x^{-n}\,dx =
   t\log\bigg(\frac{a-t}{b-t}\bigg).
\]
The quadrature formula~\eqref{cav} follows on expanding the
logarithms in powers of $t$ and comparing coefficients.

\end{document}